\documentclass{amsart}
\usepackage{latexsym}
\usepackage{amssymb}
\usepackage{euscript}

\theoremstyle{definition}

\theoremstyle{remark}

\numberwithin{equation}{section}

\begin{document}

\title{One noncommutative differential calculus coming from the inner derivation}

\author{Bo Zhao}
\address{Department of Mathematics, Washington University in St. Louis, St. Louis,MO 63130}
\curraddr{Department of Mathematics, Campus Box 1146, Washington
University in St. Louis, St. Louis, MO 63130}
\email{bzhao@math.wustl.edu}
\date{\today}

%\dedicatory{This paper is dedicated to our authors.}


\begin{abstract}
We define a noncommutative differential calculus constructed from
the inner derivation, then several relevant examples are showed.
It is of interest to note that for certain $C^*$-algebra, this
calculus is closely related to the classical one when the algebra
associates a deformation parameter.
\end{abstract}

\maketitle
\section*{Introduction}
Noncommutative geometry is developed by Alain Connes. In his work
{~\cite{AC1}}, the cyclic cohomology, which is considered to be the
classical differential calculus counterpart, is introduced. In this
paper, we are dealing with one new noncommutative differential
calculus. When considering some $C^*$-algebra $\mathcal{A}_{\theta}$
which is parametrized by a number $\theta$
($\{\mathcal{A}_{\theta}\}$ is a family of continuous field of
$C^*$-algebra, when $\theta=0$, $\mathcal{A}_{0}$ is commutative,
for instance, quantum tori, quantum plane etc), this differential
calculus will be deforming to the classical one when $\theta$ goes
to $0$. The important point to note here is that the definition is
trivial if the underlying algebra is commutative.

When combining the condition (2) of the following definition from
Rieffel {~\cite{R}}, it is very easy to see  that applying this
noncommutative differential calculus for the strict deformation
quantization algebra, it will be deforming to the classical
differential calculus.
\newtheorem{deform}{Definiton{~\cite{R}}}[section]
\begin{deform}
A strict deformation quantization of $A$ in the direction of
$\Lambda$ means an open interval $I$ of real numbers containing
$0$, together with, for each $\hbar \in I$, an associative product
$*_{\hbar}$, an involution $^{*_{\hbar}}$, and a $C^*$-norm
$||\cdot ||_\hbar$ on $A$, which for $\hbar=0$ are the original
pointwise product, complex conjugation involution, and supremum
norm, such that \begin{enumerate}
    \item for every $f\in A$, the function $\hbar\mapsto
    ||f||_{\hbar}$is continuous
    \item for every $f,g\in A$, $||(f*_{\hbar}g-g*_{\hbar}f)/i\hbar
    -\{f,g\}||_{\hbar}$ converges to $0$ as $\hbar$ goes to $0$.
\end{enumerate}
\end{deform}

The notation $\Lambda$ means a skew 2-vector field and $\{,\}$ is a
Poisson bracket. For a full treatment of the strict deformation
quantization theory, we refer the reader to {~\cite{R}}.

This paper is organized as follows. In the first section, we define
carefully this new noncommutative differential calculus and its
basic properties. In the second section, we proceed to the study of
the Dirichlet form constructed from the inner derivation. We will
show this Dirichlet form is symmetric, Markov, conservative,
completely positive and strongly local. In section 3, 4, 5, 6, 7, we
apply this calculus for $M_n(\mathbb{C})$, graph $C^*$ algebra,
quantum tori, quantum plane, and quantum Heisenberg manifolds. The
last three are the examples of the deformation quantization
$C^*$-algebra. In the case of quantum tori, the exterior derivative
is realized as a finite difference operator. It is also worth
pointing out that in section 7, we will look more closely at the
structure of the quantum Heisenberg manifolds, it turns out the von
Neumann algebra of the quantum Heisenberg manifolds is nothing but a
3-dimensional quantum tori.

I wish to thank Nik Weaver. Without his many helpful suggestions
this work could not have been done.
\section {Noncommutative differential calculus}
We are dealing with a differential calculus on some non-commutative
algebras, for instance, $C^*$-algebra, Banach algebra, von Neumann
algebra. In order to generalize this differential calculus idea, we
put the definition in the context of $C^*$-algebra. The reason to
choose $C^*$-algebra is because it always has the adjoint, and
consequently we define in the same way as from the complex
differential calculus.

Let's review first about the complex calculus. Suppose $\Omega
\subset \mathbb{C}^n$ is an open set and $f\in C^{\infty}(\Omega)$,
we have the first order calculus form
\[df=\sum_{j=1}^n \frac{\partial f}{\partial z_j}dz_j+\sum_{j=1}^n \frac{\partial f}{\partial \bar{z}_j}d\bar{z}_j,\]
here $z_j,\bar{z}_j$ are the coordinate basis. Having this complex
calculus in mind, we come up the definition below.
\newtheorem{guess}{Definiton}[section]
\begin{guess} \label{assum}
Given a noncommutative $C^*$-algebra $\mathcal{A}$, the elements
$U_1,\ldots,U_n \in \mathcal{A}$ is called a \textit{differential
basis} if the following two conditions are satisfied:
\renewcommand{\labelenumi}{(\Roman{enumi})}
\begin{enumerate}
    \item Each $U_j$ is non self-adjoint;
    \item The elements $U_1, U^{*}_1,\ldots,U_n, U^{*}_n$ mutually
    commute.
\end{enumerate}
\end{guess}
We emphasize that the condition (I) above is not necessary, at least
for certain algebra which has no adjoin, even in the $C^*$ algebra
situation. The reason we put here is to make it similar as
$z$,$\bar{z}$ from the complex calculus. The condition (II) is
crucial in order to prove $\delta^2=0$ in Proposition \ref
{bracket}. The differential basis can always be achieved from
$\mathcal{A}$, one good method is to choose  non self-adjoint
elements from any maximal abelian subalgebra of $\mathcal{A}$.
Besides the differential basis defined above, we also need to have a
derivation operator associated with the differential calculus.
\newtheorem{der}[guess]{Definiton}
\begin{der} A \textit{derivation} of a $C^*$-algebra
$\mathcal{A}$ is a linear map $\delta: \mathcal{A}\rightarrow
\mathcal{A}$ such that it satisfies the Leibniz's rule, i.e.
$\delta(xy)=\delta(x)y+x\delta(y)$. The derivation is called
\textit{inner}, if $\exists U\in \mathcal{A}$, such that $\delta
a=[U,a]=Ua-aU, \forall a\in \mathcal{A}$.
\end{der}

We find inner derivation is a good choice. The advantage of using
it lies in the fact that in certain algebra, for instance, the von
Neumann algebra, the derivation is always inner\cite{GP}. Equipped
with the differential basis and the inner derivation, we define
the first order calculus form $\delta_{U_1,\ldots,U_n}$ (we omit
the subscript $_{U_1,\ldots,U_n}$ if it's clear in the context) on
$\mathcal{A}$ as below
\[ \delta a=\sum^n_{j=1}[U_j,a]d {U_j}+\sum^n_{j=1}[
U^{*}_j,a] dU^{*}_j.
\]
The symbol $d U_j, d U^*_k$ above is merely for the notation
convenience, it can be viewed as the basis of $\delta a$ which
belongs to a direct sum of $2n$ copies of $\mathcal{A}$. For the
high order, assume $e_{\sigma(1)}\wedge \ldots \wedge
e_{\sigma(r)}=sgn(\sigma)e_1\wedge\ldots\wedge e_r$, where $\sigma$
is any permutation of $r$ numbers and $e_1,\ldots,e_r \in
\{dU_1,\ldots,dU_n, dU_1^*,\ldots,dU_n^*\}$.
 Then we put \[ \Omega
^{r}\mathcal{A}=\bigoplus_{p+q=r}\Omega^{p,q}\mathcal{A},\] where
$\Omega^{p,q}\mathcal{A}=\{\textmd {space spanned by }
dU_{i_1}\wedge\ldots\wedge dU_{i_p}\wedge
dU^{*}_{j_1}\wedge\ldots\wedge dU^{*}_{j_q} \}$.
$\Omega^{p,q}\mathcal{A}$ is called \textit{type $(p,q)$-form} and
its dimension is
\[dim\Omega^{p,q}(\mathcal{A})={2n \choose
p+q}=\frac{(2n)!}{(p+q)!(2n-p-q)!},0\leq p+q\leq 2n.\] Having
disposed of the above information, the definition for the high order
differential calculus is immediately obtained, namely if
\[\alpha= \sum a_{I,J}dU_{I} \wedge d U^{*}_{J}\in
\Omega^{p,q}\mathcal{A},\] where $I=(i_1,\ldots,i_p)$,
$J=(j_1,\ldots,j_q)$ and $dU_{I}=dU_{i_1}\wedge\ldots\wedge
dU_{i_p}$, $dU^{*}_{J}=dU^*_{j_1}\wedge\ldots\wedge dU^*_{j_q}$,
then
\[\delta\alpha=\sum_{I,J} \delta a_{I,J}\wedge dU_{I} \wedge d
U^{*}_{J}=\sum_{I,J}(\sum^n_{j=1}[U_j,
a_{I,J}]dU_j+\sum^n_{j=1}[U^{*}_j, a_{I,J}]dU^{*}_j)\wedge
dU_{I}\wedge dU^{*}_{J}.\]
\newtheorem{epro}[guess]{Definiton}
\begin{epro}
The \textit{exterior product}, or \textit{wedge product}, of a type
$(p,q)$-form $\alpha= \sum a_{I,J}dU_{I} \wedge d U^{*}_{J}$ and a
type $(r,s)$-form $\beta= \sum b_{M,N}dU_{M} \wedge d U^{*}_{N}$ is
a type $(p+s,q+t)$-form defined by
\[\alpha\wedge \beta=\sum a_{I,J}b_{M,N} dU_{I} \wedge d U^{*}_{J}\wedge dU_{M} \wedge d
U^{*}_{N}.\]
\end{epro}
The preceding definition for wedge product shows the relation
\[\alpha\wedge\beta=(-1)^{(p+q)(t+s)}\beta\wedge\alpha\] need not hold
in general if the algebra $\mathcal{A}$ is not commutative.

\newtheorem{pro}[guess]{Proposition}
\begin{pro} \label {bracket}
$ \delta^2=0 $
\end{pro}
\begin{proof}It is based on the condition (II) from definition (\ref{assum}), and we have
\[[U_j,[U_k,a]]=[U_k,[U_j,a]],[U_j^*,[U_k,a]]=[U_k,[U_j^*,a]],[U^*_j,[U^*_k,a]]=[U^*_k,[U_j^*,a]]\]. \end{proof}
\newtheorem{po1}[guess]{Propostion}
\begin{po1}
Let $\alpha$ be a $(p,q)$-form and $\beta$ be a $(s,t)$-form, then
$\delta(\alpha\wedge\beta)=\delta\alpha\wedge\beta+(-1)^{p+q}\alpha\wedge\delta\beta$.
\end{po1}
\begin{proof}This follows from the derivation property $[U,ab]=[U,a]b+a[U,b]$.
\end{proof}
For each $(p,q)$-form $\alpha= \sum a_{I,J}dU_{I} \wedge d
U^{*}_{J}$, we have a $*$ operation defined by $\alpha^*= \sum
a^*_{I,J}dU^*_{I} \wedge d U_{J}$.
\newtheorem{bim}[guess]{Proposition}
\begin{bim}\label{proj}
The set $\Omega ^{p,q}\mathcal{A}$ of type $(p,q)$-forms on
$\mathcal{A}$ is a (left)$\mathcal{A}-module$, it has the
following properties:
\begin{enumerate}
    \item If $\alpha\in \Omega^{p,q}\mathcal{A}$, then $\alpha^{*}\in
\Omega^{q,p}\mathcal{A}$.
    \item If
             $\alpha\in\Omega^{p,q}\mathcal{A}$ and $\beta\in\Omega^{r,s}\mathcal{A}$,
             then
             $\alpha\wedge\beta\in\Omega^{p+r,q+s}\mathcal{A}$.
    \item If $\alpha \in \Omega^{p,q}\mathcal{A}$, then $\delta\alpha
\in \Omega^{p+1,q}\mathcal{A}\bigoplus\Omega^{p,q+1}\mathcal{A}$.
    \item $\Omega^{p,q}\mathcal{A}=0$, if $\min\{p,q\}>n$.
\end{enumerate}
\end{bim}
\begin{proof}
The proof is straightforward by the easy computation.
\end{proof}
The point of the proceeding proposition (3) is that it allows one to
define the operators
\[\partial: \Omega^{p,q}\mathcal{A}\rightarrow
\Omega^{p+1,q}\mathcal{A}\] and
\[\partial^{*}:\Omega^{p,q}\mathcal{A}\rightarrow\Omega^{p,q+1}\mathcal{A}\]
through the direct sum projection, namely, \[\partial
\alpha=\sum_{I,J}\sum^n_{j=1}[U_j, a_{I,J}]dU_j\wedge dU_{I}\wedge
dU^{*}_{J},\]
\[\partial^{*}
\alpha=\sum_{I,J}\sum^n_{j=1}[U^*_j, a_{I,J}]dU^{*}_j\wedge
dU_{I}\wedge dU^{*}_{J}.\] From the relation
\[\delta=\partial+\partial^{*},\]  and the property
$\delta^2=0$, we have $0=(\partial+\partial^{*})^2=\partial
^2+\partial
\partial^{*}+\partial ^{*}\partial+\partial^{*2}$, hence
$\partial^{2}=\partial^{*2}=\partial \partial^{*}+\partial
^{*}\partial=0$.
\newtheorem{de1}[guess]{Definition}
\begin{de1}
A form $\alpha$ is called \emph{closed} if $\delta\alpha=0$
 and is called \emph{$\partial^*$ closed} if
$\partial^{*}\alpha=0$.
\end{de1}
Let $\mathbb{C}_{p,q}(\mathcal{A})=\{\alpha\in
\Omega^{p,q}\mathcal{A}: \delta\alpha=0\}$,
$\mathbb{HC}_{p,q}(\mathcal{A})=\{\alpha\in \Omega^{p,q}\mathcal{A}:
\partial^{*}\alpha=0\}$, then $\mathbb{C}_{p,q}(\mathcal{A})$ is the set of closed $(p,q)$-forms and
$\mathbb{HC}_{p,q}(\mathcal{A})$ is the set of $\partial^*$ closed
$(p,q)$-forms. The property $\delta ^2=0$ gives us a deRham complex,
\[0{\longrightarrow}\Omega^0\mathcal{A}\stackrel{\delta_0}{\longrightarrow}\Omega^1 \mathcal{A}
\stackrel{\delta_1}{\longrightarrow}\Omega^2 \mathcal{A}
\stackrel{\delta_2}{\longrightarrow}\ldots
\stackrel{\delta_{2n-1}}{\longrightarrow} \Omega^{2n} \mathcal{A}
\stackrel{\delta_{2n}}{\longrightarrow} 0.\] And the $k$-th
\emph{deRham cohomology} is defined by
\[H_{dR}^k(\mathcal{A})=\textrm{ker } \delta_k/\textrm{img } \delta_{k-1}.\]

The property $\partial ^{*2}=0$ gives us a  Dolbeault complex,
\[0{\longrightarrow}\Omega^{p,0}\mathcal{A}\stackrel{\partial^*_0}{\longrightarrow}\Omega^{p,1} \mathcal{A}
\stackrel{\partial^*_1}{\longrightarrow}\Omega^{p,2} \mathcal{A}
\stackrel{\partial^*_2}{\longrightarrow}\ldots
\stackrel{\partial^*_{n-1}}{\longrightarrow} \Omega^{p,n}
\mathcal{A} \stackrel{\partial^*_n}{\longrightarrow} 0.\] And the
$(p,q)$-th \emph{Delbeault cohomology} is defined by
\[H_{D}^{p,q}(\mathcal{A})=\textrm{ker } \partial^*_{q}/\textrm{img } \partial^*_{q-1}.\]

\newtheorem{dr}[guess]{Remark}
\begin{dr}
The cohomology groups depend on the choice of the differential
basis. In general, $H_{dR}^{0}(\mathcal{A})$ contains the algebra
generated by the basis. If the basis generates a maximal abelian
subalgebra, then $H_{dR}^{0}(\mathcal{A})$ is the same algebra
generated by this basis. If $A$ is a chosen differential basis and
$B\subsetneq A$, then the cohomology groups computed by $A$ and the
cohomology groups computed by $B$ need not have the subgroup
relation. This can be shown by computing $H_{dR}^1(\mathcal{A})$.
\end{dr}

\newtheorem{com}[guess]{Theorem}
\begin{com}
$H_{dR}^{0}(\mathcal{A})=H_{D}^{p,0}(\mathcal{A})=\mathbb{C}_{0,0}(\mathcal{A})$.
\end{com}
\begin{proof}
The proof is based on Fuglede-Putnam theorem which says that
$aU_i=U_i a \Leftrightarrow aU_i^*=U_i^* a$ if $U$ is normal. Then
it is a simple matter to check that
$H_{dR}^{0}(\mathcal{A})=\mathbb{C}_{0,0}(\mathcal{A})$,
$H_{D}^{p,0}(\mathcal{A})=\mathbb{HC}_{p,0}(\mathcal{A})=\mathbb{C}_{p,0}(\mathcal{A})$
and $\mathbb{HC}_{p,0}(\mathcal{A})=\mathbb{HC}_{0,0}(\mathcal{A})$,
and the proof is complete.
\end{proof}
The interest of the above theorem is that
$\mathbb{C}_{0,0}=\mathbb{HC}_{0,0}$. Considering the classical
complex differential calculus, if $f\in \mathbb{HC}_{0,0}\Rightarrow
\frac{\partial f}{\partial \bar{z}_i}=0\Rightarrow f$ is
holomorphic. So $\mathbb{HC}_{0,0}$ is the set of holomorphic
functions, and $\mathbb{C}_{0,0}$ is the set of constant functions,
therefore in the classical case $\mathbb{C}_{0,0}\neq
\mathbb{HC}_{0,0}$. The reason for this difference lies in the fact
that in the Fuglede-Putnam theorem above, $U$ need to be bounded. On
the other hand, if we require $f$ to be bounded, then by the
Liouville's theorem, $f$ is constant.

\section {Dirichlet forms} \label{diri}
Dirichlet form, in $C^*$-algebras setting, is introduced by
Albeverio and H\o egh-Krohn\cite{SR}, it shares a flavor of geometry
in the sense of Connes' noncommutative geometry\cite{AC}. For a
recent account of the theory, we refer the reader
to\cite{CF}\cite{CS}. From the classical complex case, we have the
Laplace operator defined by $\sum \frac{\partial^2}{\partial
z_i\partial \bar{z}_i}$. When transformed to the inner derivation
$[U, \cdot]$, the Laplace operator of the counterpart is
\begin{equation}\triangle=\sum^{n}_{j=1}[ U^{*}_{j}, [ U_j,\cdot]]. \label{delt}
\end{equation}

If on this  $C^*$-algebra $\mathcal{A}$, it has the following
additional assumption:
\begin{enumerate}
    \item[(III)] It has a lower semicontinuous faithful trace $\tau$.
\end{enumerate}
Then we denote by $L^2(\mathcal{A},\tau)$
($\langle,\rangle_{L^2(\mathcal{A},\tau)}$) the Hilbert space of the
GNS representation $\pi_{\tau}$ associated to $\tau$, and by
$L^\infty(\mathcal{A},\tau)$ or $\mathcal{M}$ the von Neumann
algebra $\pi_{\tau}(\mathcal{A})''$ in
$\mathcal{B}(L^2(\mathcal{A},\tau))$ generated by $\mathcal{A}$ in
the GNS representation. $1_{\mathcal{M}}$ stands for the unit of
$\mathcal{M}$.

\newtheorem{de2}[guess]{Definition}
\begin{de2}
Given a strongly continuous semigroup $\Phi_t (t\in R^+)$ of
operators defined on $L^\infty(\mathcal{A},\tau)$,
\begin{enumerate}
    \item it is
\emph{symmetric}, if $\tau( \Phi_t(x)y)=\tau(x \Phi_t(y))$.
    \item it is \emph{Markov}, if $0\leq x\leq 1_{\mathcal{M}}$ implies that $0\leq
\Phi_t(x)\leq 1_{\mathcal{M}}$.
    \item it is \emph{conservative}, if
$\Phi_t(1_{\mathcal{M}})=1_{\mathcal{M}}$.
    \item it is \emph{completely positive}, if for any $n$ we have $\sum_{i,j=1}^n
b_i^* \Phi_t(a_i^* a_j)b_j\geq 0$ where $a_i,b_i\in \mathcal{M}
,i=1,\ldots,n$.
\end{enumerate}
\end{de2}
\newtheorem{pr2}[guess]{Proposition}
\begin{pr2} \label{cp}
$\triangle$ is the generator of a norm continuous, completely
positive, symmetric, conservative Markov semigroup.
\end{pr2}
\begin{proof}Let $\Phi_t=e^{-t\triangle}$, then $\Phi_t$ is norm
continuous because $\triangle$ is bounded.\\
$\Phi_t$ is symmetric and conservative:\\
$\tau(\triangle(a)b)=\sum_{j=1}^n\tau(
[U_j^{*},[U_j,a]]b)=\sum_{j=1}^n\tau(a [U_j^{*},[U_j,b]])=\tau(a\triangle(b))\Rightarrow \tau(\triangle^n(a)b)=\tau(a\triangle^n(b))\Rightarrow\tau(\Phi_t(a)b)=\tau(a\Phi_t(b))$, also $\Phi_t(1_{\mathcal{M}})=1_{\mathcal{M}}$.\\
$\Phi_t$ is completely positive:\\
$-\triangle a=\sum_{j=1}^n(U_j^*aU_j+U_j a
U_j^*-U_j^*U_ja-aU_jU_j^*)$, let $K_1(a)=\sum_{j=1}^n(U_j^*aU_j+U_j
a U_j^*)$ and $K_2(a)=Aa+aA$ where $A=-\sum_{j=1}^nU_j^*U_j=A^*$. As
$K_1$ is completely positive, so is $K_1^n$, hence $e^{tK_1}$ is
completely positive. $(1+\frac{tx}{n})^n\rightarrow e^{tx}$, for
$x\in\mathbb{R}$. From the functional calculus theorem,
$A_m=(1+\frac{tK_2}{m})^m\rightarrow e^{tK_2}$ in norm. Since
$\sum_{i,j=1}^n b_i^* K_2(a_i^* a_j)b_j$ is self-adjoint,
$1+\frac{tK_2}{m}$ is completely positive when $m$ large enough$
\Rightarrow A_m$ is completely positive$\Rightarrow e^{tK_2}$ is
completely positive. Then from the Trotter product formula,
$(e^{\frac{t}{n}K_1}e^{\frac{t}{n}K_2})^n\rightarrow
e^{-t\triangle}$ in the strong operator topology. Hence $\Phi_t$ is
completely positive.\\
$\Phi_t$ is Markov:\\
From above, in particular, $\Phi_t$ is positive, then for $0\leq
a\leq 1_{\mathcal{M}}$, $e^{-t\triangle}(a-1_{\mathcal{M}})\leq
0\Rightarrow 0\leq e^{-t\triangle}a \leq 1_{\mathcal{M}}$.
\end{proof}

From ~\cite{SR}, $\mathcal{E}(a,b)=\langle\triangle^{\frac{1}{2}} a,
\triangle^{\frac{1}{2}} b\rangle_{L^2(\mathcal{A},\tau)}$ is the
Dirichlet form with the corresponding generator $\triangle$. With
this Dirichlet form, we can construct a $C^*$-Hilbert
$\mathcal{A}$-bimodule $\mathcal{A}\otimes \mathcal{A}$~\cite{S} and
a sesquilinear form with values in $\mathcal{A}$ by the formula
\begin{equation}
\langle a\otimes b, c\otimes d
\rangle_{\mathcal{A}}=b^*(a^*\triangle
c+\triangle(a^*)c-\triangle(a^* c))d. \label{loca}
\end{equation}

\newtheorem{de3}[guess]{Definition{\cite{S}}}
\begin{de3}
A completely positive Markov semigroup $(\Phi_t)_{t\geq 0}$ and
its infinitesimal generator $\triangle$ are \emph{strongly local}
if there exists a Hilbert space $\mathbb{H}$ and an isometry $W$
in $\mathcal{L}_{\mathcal{A}}(\mathcal{A}\otimes
\mathcal{A},\mathbb{H}\otimes\mathcal{A})$.
\end{de3}
Replacing $\triangle$ defined in (\ref{delt}) to (\ref{loca}), we
have
\[\langle a\otimes b, c\otimes d \rangle_{\mathcal{A}}=\sum_{j=1}^n b^*([U_j,a]^* [U_j,c]+
[U^*_j,a]^*[U^*_j,c] )d.\] Take $\mathbb{H}$ to be
$\mathbb{C}^{2n}$, and define the mapping $W: \mathcal{A}\otimes
\mathcal{A} \rightarrow \mathbb{C}^{2n}\otimes\mathcal{A}$ by
\[a\otimes b \mapsto [U_1,a] b\oplus\ldots
\oplus [U_n^*,a]b,\] and the $\mathcal{A}$-bimodule structure by
\[\langle [ U_1,a] b\oplus\ldots \oplus [U_n^*,a] b , [U_1,c]
d\oplus\ldots \oplus [U_n^*,c]d \rangle_{\mathcal{A}}=b^*
\sum_{j=1}^n ([U_j,a]^*[U_j,c]+ [U^*_j,a]^*[U^*_j,c])d.\] It's not
hard to see $W$ turns out to be an isometry, and this is precisely
the assertion of the following theorem.

\newtheorem{sll}[guess]{Theorem}
\begin{sll}
$\triangle$ defined in (\ref{delt}) is strongly local.
\end{sll}

Since $\Phi_t$ is conservative, by theorem (4.7)\cite{CS}, this
Dirichlet form has a representation
\begin{equation}
\mathcal{E}(a,a)=\tau(\langle \delta a, \delta a
\rangle_{\mathcal{A}}).  \label{di}
\end{equation}

\section {Matrix algebra} \label{matrix}
For the $n \times n$ matrix algebra $M_n(\mathbb{C})$, the diagonal
matrices is its maximal abelian subalgebra. The elements
$p_1,\ldots,p_n$ where $p_j$ is the $n\times n$ matrix with $1$ in
entry $(j,j)$ and $0$ elsewhere are the generators of the diagonal
matrices. Therefore $p_1,\ldots,p_n$ can be chosen as the
differential basis. Thus we actually have constructed the
noncommutative differential calculus on $M_n(\mathbb{C})$. Its first
order differential calculus becomes
\[\delta_n a =[p_1,a]dp_1+\ldots+[p_n,a]dp_n.\] The inner derivation
has the form
\[[p_j,a]=\begin{pmatrix} & && -a_{1,j} & & &   \\
                  & && \vdots & & &   \\
                  & && -a_{j-1,j} & & &   \\
                  a_{j,1}&\ldots &a_{j,j-1}& 0 &a_{j,j+1} &\ldots &a_{j,n}   \\
                  & && -a_{j+1,j} & & &   \\
                  & && \vdots & & &   \\
                  & && -a_{n,j} & & &   \\
                    \end{pmatrix}.\]

For the matrix algebra $M_n(\mathbb{C})$, we have the matrix trace
$Tr$ defined by $Tr(a)=\sum_{i=1}^n a_{i,i}$. Since $Tr$ satisfies
the condition (III) from section ~\ref{diri}, we actually have
defined a Dirichlet form on $M_n(\mathbb{C})$. Now let's compute
it explicitly. In the proposition below, the notation $\tilde{a}$
 means $\tilde{a}_{i,j}=a_{i,j}$ when $i\neq j$ and $\tilde
{a}_{j,j}=0$, i.e., $\tilde{a}$ kills the diagonal entries of $a$.
We will keep using this notation in this section.
\newtheorem{di}[guess]{Proposition}
\begin{di}
The Dirichlet form on $M_n(\mathbb{C})$ constructed from $Tr$ and
the Laplace operator $\sum_{j=1}^n [p_j,[p_j,\cdot]]$ is
\[\mathcal{E}_{n}(a,b)=2 Tr(\tilde{a}^{*} \tilde{b}).\]
\end{di}
\begin{proof}
Given the Laplace operator and $Tr$, from (\ref{di}), this Dirichlet
form defined on $M_n(\mathbb{C})$ is
\[\mathcal{E}_{n}(a,b)=Tr(\sum_{j=1}^n[p_j,a]^*[p_j,b]).\]
As \[Tr([p_j,a]^*[p_j,b])=\sum_{i=1,i\neq j}^n
(a^*_{i,j}b_{i,j}+a^*_{j,i}b_{j,i})=\sum_{i=1}^n
(\tilde{a}^*_{i,j}\tilde{b}_{i,j}+\tilde{a}^*_{j,i}\tilde{b}_{j,i}),\]
therefore, $\mathcal{E}_{n}(a,b)=
\sum_{i,j}(\tilde{a}^*_{i,j}\tilde{b}_{i,j}+\tilde{a}^*_{j,i}\tilde{b}_{j,i})=2
Tr(\tilde{a}^*\tilde{b}).$
\end{proof}

So far, we have defined the Dirichlet form $\mathcal{E}_{n}$ on
$M_n(\mathbb{C})$. Regard $M_n(\mathbb{C})$ as a subalgebra of
$M_{n+1}(\mathbb{C})$ via the embedding
\[a\rightarrow \begin{pmatrix} a & 0   \\
                  0 & 0   \\
                    \end{pmatrix},\]
we see that $\delta_{n+1}$ is compatible with $\delta_{n}$, i.e.,
$\delta_{n+1}(a)=\delta_n(a)$, $\forall a\in M_n(\mathbb{C})$, and
$\mathcal{E}_{n+1}$ is compatible with $\mathcal{E}_{n}$, i.e.,
$\mathcal{E}_{n+1}(a,b)=\mathcal{E}_{n}(a,b), \forall a,b\in
M_n(\mathbb{C})$. We are thus led to define the differential
calculus, Dirichlet form on $\cup_{n} M_n(\mathbb{C})$. As the norm
closure of $\cup_{n} M_n(\mathbb{C})$ is $\mathcal{K}(l^2)$, i.e,
the compact operators on $l^2$, we have the following theorem after
this extension. In the theorem below, we use the notation
$\mathcal{HS}(l^2)$ to represent the Hilbert-Schmidt operators on
$l^2$, i.e., $x\in \mathcal{HS}(l^2)\Leftrightarrow
Tr(x^*x)<\infty$.

\newtheorem{ex}[guess]{Theorem}
\begin{ex}\label{dir}
$\mathcal{E}(a,b)=2Tr(\tilde{a}^*\tilde{b})$ is a Dirichlet form on
$\mathcal{HS}(l^2)$, the corresponding Laplace operator is
$\sum_{j=1}^{\infty} [p_j,[p_j,\cdot]]$. The restriction of
$\mathcal{E}$ on $M_n(\mathbb{C})$ is $\mathcal{E}_n$.
\end{ex}
\begin{proof}
We see at once that
 $L^{\infty}(\mathcal{K}(l^2),Tr)=\mathcal{B}(l^2)$,
$L^{2}(\mathcal{K}(l^2),Tr)=\mathcal{HS}(l^2)$, and
$Dom(\mathcal{E})=\mathcal{HS}(l^2)$. If we let
$\Delta=\sum_{j=1}^{\infty} [p_j,[p_j,\cdot]]$, then
$\Delta(a)=\tilde{a}$ for $a\in \mathcal{B}(l^2)$. In this way, the
semigroup $\Phi_t(a)=e^{-t\Delta}a=a-\tilde{a}+e^{-t}\tilde{a}$. An
easy computation shows $Tr(\tilde{a}^*\tilde{b})=Tr(\Delta(a)^*b)$
and
\[-\Delta(a)=2\sum_{i=1}^{\infty}p_j ap_j-2a.\] By using the same
argument as from proposition (\ref{cp}), $\triangle$ is seen to be
the generator of a symmetric, conservative, completely positive
Markov semigroup, which completes the proof.
\end{proof}

If we regard $M_{2^n}(\mathbb{C})$ as a subalgebra of $M_{2^{n+1}}(\mathbb{C})$  via the embedding \[a\rightarrow \begin{pmatrix} a & 0   \\
                  0 & a   \\
                    \end{pmatrix},\]
and consider the fact there exists a unique tracial state $tr$, we
can perform the same construction as above and define a Dirichlet
form on hyperfinite $II_1$ factor. ~\cite{bz} provides a detailed
study for this construction.

The following is the deRham complex of $M_n(\mathbb{C})$ with the
differential basis $p_1,\ldots,p_n$.
\[0\stackrel{\delta}{\longrightarrow}M_n(\mathbb{C})\stackrel{\delta}{\longrightarrow}\Omega^1 M_n(\mathbb{C}) \stackrel{\delta}{\longrightarrow}\Omega^2 M_n(\mathbb{C})
\stackrel{\delta}{\longrightarrow}\ldots
\stackrel{\delta}{\longrightarrow} \Omega^n M_n(\mathbb{C})
\stackrel{\delta}{\longrightarrow} 0.\] The $0$-th deRham cohomology
group of $M_n(\mathbb{C})$, which we will show from section
~\ref{graph}, is
\[ H^0_{dR}(M_n(\mathbb{C}))=\mathbb{C}^n.\]

\section {Graph $C^*$-algebras} ~\label{graph}
For a more detailed introduction to graph $C^*$-algebras we refer
to ~\cite{TDIW},\cite{ADI} and the reference therein. A directed
graph $E=(E^0,E^1,r,s)$ consists of countable sets $E^0$ of
vertices and $E^1$ of edges, and maps $r,s:E^1\rightarrow E^0$
identifying the range and source of each edge. The graph is
\textit{row-finite} if each vertex emits at most finitely many
edges. We write $E^n$ for the set of paths
$\mu=\mu_1\mu_2\ldots\mu_n$ of length $|\mu|:=n$; that is,
sequences of edges $\mu_i$ such that $r(\mu_i)=s(\mu_{i+1})$ for
$1\leq i<n$. The map $r,s$ extend to $E^*:=\cup_{n\geq0}E^n$ in an
obvious way, and $s$ extends to the set $E^{\infty}$ of infinite
paths $\mu=\mu_1\mu_2\ldots$. A \textit{sink} is a vertex $v\in
E^0$ with $s^{-1}(v)=\emptyset$, a \textit{source} is a vertex
$w\in E^0$ with $r^{-1}(w)=\emptyset$.

A \textit{Cuntz-Krieger E-family} in a $C^*$-algebra $B$ consists
mutually orthogonal projections $\{p_v:v\in E^0\}$ and partial
isometries $\{s_e: e\in E^1\}$ satisfying the Cuntz-Krieger
relations
\[s^*_e s_e=p_{r(e)} \textrm{ for }e\in E^1\textrm{ and }p_v=\sum_{\{e:s(e)=v\}}s_es^*_e\textrm{ whenever }v\textrm{ is not a sink}.\]
It is proved in ~\cite{ADI} that there is a universal $C^*$-algebra
$C^*(E)$ generated by a non-zero Cuntz-Krieger E-family
$\{s_e,p_v\}$. A product $s_{\mu}:=s_{\mu_1}s_{\mu_2}\ldots
s_{\mu_n}$ is non-zero precisely when $\mu=\mu_1\mu_2\ldots\mu_n$ is
a path in $E^n$. Since the Cuntz-Krieger relations imply that the
projections $s_es_e^*$ are also mutually orthogonal, we have
$s^*_es_f=0$ unless $e=f$, and words in $\{s_e,s_f^*\}$ collapse to
products of the form $s_{\mu}s^*_{\nu}$ for $\mu,\nu\in E^*$
satisfying $r(\mu)=r(\nu)$. Indeed, because the family
$\{s_{\mu}s^*_{\nu}\}$ is closed under multiplication and
involution, we have
\begin{equation}
C^*(E)=\overline{\textrm{span}}\{s_{\mu}s^*_{\nu}: \mu,\nu \in
E^*\textrm{ and }r(\mu)=r(\nu)\}.\label{dense}
\end{equation}
We adopt the conventions that vertices are paths of length $0$,
that $s_{v}:=p_v$ for $v\in E^0$, and all the paths $\mu,\nu$
appearing in (\ref{dense}) are non-empty; we recover $s_{\mu}$,
for example, by taking $\nu=r(\mu)$, so that
$s_{\mu}s^*_{\nu}=s_{\mu}p_{r(\mu)}=s_{\mu}$.

For simplicity, we will only consider the graph which has the finite
vertices. It is then easily seen that $\{p_v: v\in E^0\}$ satisfies
the condition (I)(II) from definition ~\ref{assum}, in this way, we
have defined the noncommutative differential calculus on $C^*(E)$.
\newtheorem{rel}[guess]{Lemma}
\begin{rel}
Fix a vertex $v_0\in E^0$ and a path $\mu$. If $|\mu|=n (n\neq 0)$,
let $\mu=\mu_1\mu_2\ldots\mu_n$, otherwise  $\mu$ is a vertex.
\begin{enumerate}
    \item
if $s(\mu) \neq r(\mu)$, then
\begin{equation*}
[p_{v_0}, s_{\mu}]=
\begin{cases} s_{\mu} & \text{if $v_0=s(\mu)$,}\\
-s_{\mu} & \text{if $v_0=r(\mu)$,}\\
0 &\text{otherwise.}
\end{cases}
\end{equation*}

    \item if $s(\mu)=r(\mu)$, i.e., path $\mu$ is a loop,
    then
$[p_{v_0}, s_{\mu}]=0$.
 \item $\delta (s_{\mu})= s_{\mu}dp_{s(\mu)}-
s_{\mu}dp_{r(\mu)}$.
\end{enumerate}
\end{rel}
\begin{proof}
The case $|\mu|=0$ is easy. So we assume $|\mu|=n$, and
$\mu=\mu_1\mu_2\ldots\mu_n$. From the identification that
$p_{v_0}=s_{v_0}$, we know $s_{v_{0}}s_{\mu}=s_{\mu}$ iff
$s(\mu)=v_0$ and $s_{\mu}s_{v_{0}}=s_{\mu}$ iff $r(\mu)=v_0$, which
proves the lemma.
\end{proof}
The remainder of this section will be devoted to discuss the $0$-th
deRham cohomology of $C^*(E)$.
\newtheorem{h0}[guess]{Proposition}
\begin{h0}
$\mathbb{C}_{0,0}=\{s_{\mu}s^*_{\nu}: \mu,\nu\in E^*, r(\mu)=r(\nu),
s(\mu)=s(\nu))\}$
\end{h0}
\begin{proof}
The above lemma and relation $[p_{v_0}, s_{\mu}^*]=-[p_{v_0},
s_{\mu}]^*$ shows
\begin{eqnarray*}
          \delta(s_{\mu}s^*_{\nu})&=& \delta(s_{\mu})s_{\nu}^*+s_{\mu}\delta(s^*_{\nu})
          \\& = &
            s_{\mu}s^*_{\nu}dp_{s(\mu)}-
            s_{\mu}s^*_{\nu}dp_{r(\mu)}-s_{\mu}(s^*_{\nu}dp_{s(\nu)}- s^*_{\nu}dp_{r(\nu)})
                \\ & = &
            s_{\mu}s^*_{\nu}dp_{s(\mu)}
            -s_{\mu}s^*_{\nu}dp_{s(\nu)}
                \end{eqnarray*}
\end{proof}

\newtheorem{h1}[guess]{Corollary}
\begin{h1}
If $E$ has no loops, then $\mathbb{C}_{0,0}=\{s_{\mu}s^*_{\mu}:
\mu\in E^*\}$.
\end{h1}
Notice, $s_{\mu}$ is a nonzero partial isometry with
$s_{\mu}s^*_{\mu}\leq p_{s(\mu)}$. Given a path $\mu$, the next
proposition gives a criteria for when $s_{\mu}s^*_{\mu}= p_{s(\mu)}$
is true.

\newtheorem{path}[guess]{Proposition}
\begin{path}
Suppose $\mu$ is a path and $\mu=\mu_1\mu_2\ldots\mu_n$, then
$s_{\mu}s^*_{\mu}=p_{s(\mu)}$ iff $\{e\in E^1:
s(e)=s(\mu_j)\}=\{\mu_j\}$ for $j=1,\ldots, n-1$. Intuitively, it is
the kind of path, on which all the nodes, except the range node,
have only one exit.
\end{path}
\begin{proof}
The proof follows from the Cuntz-Krieger relations
\[s^*_e s_e=p_{r(e)} \textrm{ for }e\in E^1\textrm{ and }p_v=\sum_{\{e:s(e)=v\}}s_es^*_e\textrm{ whenever }v\textrm{ is not a sink}.\]
\end{proof}
The above two propositions provide a way to get the deRham
cohomology $H_{dR}^0(C^*(E))$, which has special meaning behind the
graph itself. If the directed graph $E$ has no loops, then
$H^0_{dR}(C^*(E))\cong \mathbb{C}^{m}$, where $m\geq \#
\textrm{vertices}$. In particular, if the graph $E$ is a tree with
only one sink, then $m= \#\textrm{vertices}$. The matrix algebra
$M_n(\mathbb{C})$ discussed in section ~\ref{matrix} is a tree with
one root and $n-1$ nodes connecting to the root. And we have
$H_{dR}^0(M_n(\mathbb{C})))=n$. When $E$ has loops,
$H_{dR}^0(C^*(E))$ is a bit more complicated, say if $E$ contains a
simple loop with the property that this loop doesn't have exit, then
$H^0_{dR}(C^*(E))$ contains $C(\mathbb{T})$. (the loop is simple
when its vertices are distinct)

\section {2n-dimensional quantum tori} \label{tori}
Let's look at the $2n$-dimensional quantum tori, $U_1,\ldots,
U_{2n}$ are its unitary generators which subject to the relation
$U_k U_j=e^{i\theta_{j,k}}U_j U_k$. $\Theta=(\theta_{j,k})$ is a
skew symmetric matrix which by the general spectral theorem, it
can be orthogonal transformed to a block matrix,  with the entry
(we use the same letter) $\theta_{j,k}=0$ except possibly these
$\theta_{2j-1,2j}$ and $\theta_{2j,2j-1} (j=1,\ldots,n)$. Then the
unitary operators $U_{2j-1}, j={1,\ldots,n}$ satisfy the
conditions (I) (II) from definition (\ref{assum}). So we have
defined the noncommutative differential calculus on this quantum
tori. Its first order differential form is
\[\delta A=\sum_{j=1}^n [U_{2j-1}, A]dU_{2j-1}+\sum_{j=1}^n [U^*_{2j-1}, A]dU^*_{2j-1}.\]

Take for simplicity, let's assume n=1, and the unitary generators
$U,V$ satisfies the relation $UV=e^{i\theta}VU$. In order to compute
its cohomology groups, we need first to understand the set
$\mathbb{C}_{0,0}=\{A\in\mathcal{A}: [ U,A]=0\}$.

Define $\hat{\theta}_{s,t}(U^kV^l)=e^{-i(sk+tl)}U^kV^l$ $((s,t)\in
\mathbb{R}^2)$, then ~\cite{N} shows $\hat{\theta}$ defines an
action of $\mathbb{R}^2$ by automorphism of $\mathcal{A}$.
\newtheorem{th1}[guess]{Theorem}
\begin{th1}
If $\theta/\pi$ is irrational, then $\mathbb{C}_{0,0}=<U>$(the
algebra generated by $U$).
\\
If $\theta/\pi=\frac{p}{q}(p\neq2)$, then $\mathbb{C}_{0,0}=<U,V^{2q}>=\{A\in\mathcal{A}: \hat{\theta}_{0,\frac{\pi}{q}}(A)=A\}$. \\
If $\theta/\pi=\frac{2}{q}$, then
$\mathbb{C}_{0,0}=<U,V^q>=\{A\in\mathcal{A}:
\hat{\theta}_{0,\frac{2\pi}{q}}(A)=A\}$.
\end{th1}
\begin{proof} The prove is based on the concept of Fourier analysis method, which is already shown from \cite{N}. The Fourier
coefficient of $A$ is $a_{k,l}(A)=\langle
Ae_{0,0},e_{k,l}\rangle$, thereafter
\[a_{k,l}(UA)=\langle UAe_{0,0},e_{k,l}\rangle=\langle Ae_{0,0},e^{-i\theta l/2}e_{k-1,l}\rangle=e^{i\theta l/2}a_{k-1,l}(A),\]
\[a_{k,l}(AU)= \langle AUe_{0,0},e_{k,l}\rangle= \langle Ae_{1,0},e_{k,l} \rangle=e^{i\theta (-l)/2}a_{k-1,l}(A),\]
hence if $A\in \mathbb{C}_{0,0}$, then $UA=AU$, $0=e^{i\theta
l/2}(1-e^{-i\theta l})a_{k-1,l}$, and the proof is complete.
\end{proof}
The following is the deRham complex of $\mathcal{A}$ from this
differential operator $\delta$,
\[0\stackrel{\delta}{\longrightarrow}\mathcal{A}\stackrel{\delta}{\longrightarrow}\Omega^1\mathcal{A} \stackrel{\delta}{\longrightarrow}\Omega^2\mathcal{A}
\stackrel{\delta}{\longrightarrow} 0.\] From above theorem, when
$\theta$ is irrational,
\[ H^0_{dR}(\mathcal{A})=<U>\cong C(\mathbb{T}).\]

The principal significance of the above introduced differential
calculus is that it allows to deform to the classical one. The
remainder of this section will be devoted to show it.

First let's investigate the geometric meaning behind the derivation
$\delta$. Suppose an element $a\in \mathcal{A}$ for the moment has
the form $a=f(U,V)=\sum a_{m,n}U^m V^n$. In fact, it is true that
the above elements are dense in $\mathcal{A}$. Take $U$ as the
differential basis first, then the first order form is
\begin{eqnarray*}
\delta_{U} f(U,V) & = & U(f(U,V)-f(U,e^{-i\theta}V)) dU  +
U^*(f(U,V)-f(U,e^{i\theta}V)) d U^*
\\ & = & U(a-\hat{\theta}_{-\theta,0}(a)) dU  +
U^*(a-\hat{\theta}_{\theta,0}(a)) d U^*
\end{eqnarray*}
As $\hat{\theta}_{-\theta,0}$ is continuous in norm, for any $a\in
\mathcal{A}$, we actually have
\[\delta_{U}a=U(a-\hat{\theta}_{-\theta,0}(a)) dU  +
U^*(a-\hat{\theta}_{\theta,0}(a)) d U^*.\]

More general, we could take $U^{k_1}V^{k_2}$ as the differential
basis. A slight change in the above observation actually leads to
the following theorem.

\newtheorem{th2}[guess]{Theorem}
\begin{th2}
Given $a\in \mathcal{A}$, the first order form by taking
$U^{k_1}V^{k_2}$ as the differential basis is
\begin{eqnarray*}
\delta_{U^{k_1}V^{k_2}} a & = &
U^{k_1}V^{k_2}(a-\hat{\theta}_{-k_1\theta,k_2\theta}(a))dU^{k_1}V^{k_2}
\\ & + & (U^{k_1}V^{k_2})^*(a-\hat{\theta}_{k_1\theta,-k_2\theta}(a))d
(U^{k_1}V^{k_2})^*.
\end{eqnarray*}
\end{th2}

Interestingly, this derivation involves finite difference rather
than derivations. Namely, $[U^{k_1}V^{k_2}, a]$ is the changing in
$a$ along the $V$ direction by $-k_1\theta$ difference and the $U$
direction by $k_2\theta$ difference. This is a purely noncommutative
phenomenon because difference quotients can't be used as the basis
of a differential calculus in the commutative case.

Next, coming back to the $2n$-dimensional quantum tori, based on the
above theorem, we will show that the above first order form can be
deformed to the classical one on $T^n$. In the sequel, $\theta_j$
denotes $\theta_{2j-1,2j}\neq 0, j=1,\ldots,n$. We consider the
derivative operator $\partial$ and $C^*$ subalgebra $\mathcal{B}$ of
$\mathcal{A}$ generated by $U_2,U_4,\ldots,U_{2n}$. The relation
matrix $\Theta$ tells us $\mathcal{B}$ is commutative, and hence it
is $*$ isomorphic to $C(\mathbb{T}^n)$.

For our purpose, instead of taking the differential basis from the
beginning of the section, we slightly change to
 $\theta_1^{-1}U_1,\ldots,\theta_j^{-1}U_{2j-1},\ldots,\theta_n^{-1}U_{2n-1}$. Indeed, this again satisfies condition (I) and (II). For $a=f(U_2,U_4,\ldots,U_{2n})\in
 \mathcal{B}$,
\[\partial a=\sum_{j=1}^{n} [\theta_j^{-1} U_{2j-1},a] dU_{2j-1}=\sum_{j=1}^{n} \frac{1}{\theta_j}[ U_{2j-1},a]dU_{2j-1}.\]
Then, \[ \partial f= \sum_{j=1}^{n}
\frac{1}{\theta_j}(f(U_2,\ldots,U_{2j}e^{i\theta_j},\ldots,U_{2n})-f(U_2,\ldots,U_{2j},\ldots,U_{2n}))U_{2j-1}dU_{2j-1}.\]
Since \[f(U_2,U_4,\ldots,U_{2n})\cong
f(e^{ix_2},e^{ix_4},\ldots,e^{ix_{2n}}),\]
\begin{eqnarray*}
&& \partial f
\\ &\cong&
\sum_{j=1}^{n}
\frac{f(e^{ix_2},\ldots,e^{i(x_{2j}+\theta_j)},\ldots,e^{ix_{2n}})-f(e^{ix_2},\ldots,e^{ix_{2j}},\ldots,e^{ix_{2n}})}{\theta_j}U_{2j-1}dU_{2j-1}
\\& \rightarrow &
\sum_{j=1}^{n} \frac{\partial
f(e^{ix_2},\ldots,e^{ix_{2j}},\ldots,e^{ix_{2n}})}{\partial
e^{ix_{2j}}}ie^{ix_{2j}}U_{2j-1}dU_{2j-1}
\\ & = &
\sum_{j=1}^{n} \frac{\partial
f(U_2,\ldots,U_{2j},\ldots,U_{2n})}{\partial
U_{2j}}iU_{2j}U_{2j-1}dU_{2j-1}
\\ & \cong &
\sum_{j=1}^{n} \frac{\partial
f(U_2,\ldots,U_{2j},\ldots,U_{2n})}{\partial U_{2j}}dU_{2j}
\end{eqnarray*}
as $\theta_j\rightarrow 0, j=1,\ldots,n$, which turns out to be
quite similar as the classical first order differential calculus on
$\mathbb{T}^n$. The same argument works for any high order. The
example from quantum tori demonstrates rather strikingly that this
noncommutative differential calculus is closely related to the
classical one.

\section {Quantum Plane}
In physics, it has the "quantization procedure" which from a
classical Hamiltonian \[H(p_1,q_1,\ldots,p_n,q_n)\] on the
2n-dimensional phase space $\mathbb{R}^{2n}$ to the quantum
mechanical Hamiltonian version \[H(P_1,Q_1,\ldots,P_n,Q_n).\] Here
$Q_j$ denotes the multiplication operator on $L^2(\mathbb{R}^n)$
corresponding to the coordinate mapping $x_j$, and $P_j=-i\hbar
\frac{\partial}{\partial x_j}$ is a partial derivative operator.
The procedure for getting the differential calculus on quantum
plane is quite similar as from the quantum tori, we continue in
the fashion by first investigating the simple case $n=1$.
\newtheorem{qp}[guess]{Definition}
\begin{qp}
A canonical pair is informally described as a pair of self-adjoint
operators $P,Q$ on a Hilbert space $\mathcal{H}$, satisfying
(Heisenberg's) canonical commutation relation(CCR)
\[[P,Q]=-i\hbar I\]
\end{qp}
The momentum operator defined on $L^2(\mathbb{R}^2)$ by
$Pf(x,y)=yf(x,y)-i\frac{\hbar}{2}\frac{\partial f(x,y)}{\partial
x}$ and the position operator defined by
$Qf(x,y)=xf(x,y)+i\frac{\hbar}{2}\frac{\partial f(x,y)}{\partial
y}$ satisfy the CCR relation.

By using the Fourier transform, we obtain the quantum observable
\[f(P,Q)=\frac{1}{2\pi}\int\hat{f}(t_1,t_2)e^{i(t_1P+t_2Q)}dt_1dt_2,\]
where $f$ is in the Schwartz class $\mathcal{S}(\mathbb{R}^2)$.

Fix $k_1,k_2\in \mathbb{R}$, we do the differential calculus with
the basis $\frac{1}{\hbar}e^{ik_1P+ik_2Q}$. Then \[\delta
f(P,Q)=[\frac{e^{ik_1P+ik_2Q}}{\hbar},
f(P,Q)]de^{ik_1P+ik_2Q}+[\frac{e^{-ik_1P-ik_2Q}}{\hbar},
f(P,Q)]de^{-ik_1P-ik_2Q}\] where $f\in \mathcal{S}(\mathbb{R}^2)$.
\begin{eqnarray*}
& & \lim_{\hbar\rightarrow
0}[\frac{e^{ik_1P+ik_2Q}}{\hbar},f(P,Q)]
\\ & = &
\lim_{\hbar\rightarrow
0}\frac{1}{2\pi}\int\hat{f}(t_1,t_2)\frac{1}{\hbar}[e^{ik_1P+k_2iQ},e^{i(t_1P+t_2Q)}]dt_1dt_2
\\ & = &
\lim_{\hbar\rightarrow 0}\frac{1}{2\pi
\hbar}\int\hat{f}(t_1,t_2)e^{i(t_1P+t_2Q)}e^{ik_1P+ik_2Q}(e^{ik_1t_2\hbar-ik_2t_1\hbar}-1)dt_1dt_2
\\ & = &
\frac{1}{2\pi}\int\hat{f}(t_1,t_2)e^{i(t_1P+t_2Q)}i(k_1t_2-k_2t_1)dt_1dt_2e^{ik_1P+ik_2Q}
\\ & = &
 (k_1\frac{\partial f(P,Q)}{\partial Q}-k_2 \frac{\partial
f(P,Q)}{\partial P})e^{ik_2P+ik_2Q}.
\end{eqnarray*}
In conclusion, when $\hbar\rightarrow 0$, $\delta f(P,Q)$ deforms
to
\begin{eqnarray*}
& &(k_1\frac{\partial f(P,Q)}{\partial Q}-k_2 \frac{\partial
f(P,Q)}{\partial P})e^{ik_2P+ik_2Q}de^{ik_1P+ik_2Q}
\\& - &(k_1\frac{\partial f(P,Q)}{\partial Q}-k_2 \frac{\partial
f(P,Q)}{\partial P})e^{-ik_1P-ik_2Q}de^{-ik_1P-ik_2Q}.
\end{eqnarray*}
The above deformed differential calculus, comparing to the
classical one, is the linear sum in two directions.

Having disposed of the simple case $n=1$, we now return to the
general 2n-dimensional quantum plane. Quite the same way as from
the quantum tori, we take the differential basis to be
$\frac{1}{\hbar}e^{iP_1},\ldots,\frac{1}{\hbar}e^{iP_n}$, it
satisfies the condition (I) and (II) from definition
(\ref{assum}). We consider only the derivative operator $\partial$
and subalgebra $\mathcal{B}$ generated by $Q_1, Q_2,\ldots, Q_n$.
By the commutative relations of $Q_j$, $\mathcal{B}$ is isomorphic
to $C_0(\mathbb{R}^n)$. We start doing the first order calculus of
the element $f(Q_1,\ldots,Q_n)$ where $f\in
\mathcal{S}(\mathbb{R}^n)$.
\begin{eqnarray*}
& & \lim_{\hbar\rightarrow 0}
\partial_{\frac{1}{\hbar}e^{iP_1},\ldots,\frac{1}{\hbar}e^{iP_n}}
f(Q_1,\ldots,Q_n)
\\ &=& \lim_{\hbar\rightarrow 0} \sum_{j=1}^n
[\frac{1}{\hbar}e^{iP_j},f(Q_1,\ldots,Q_n)]de^{iP_j}
\\ &=&
\lim_{\hbar\rightarrow 0}
\sum_{j=1}^n\frac{1}{2\pi}\int\hat{f}(t_1,\ldots,t_n)e^{i(t_1Q_1+\ldots+t_nQ_n)}\frac{e^{it_j\hbar}-1}{\hbar}dt_1\ldots
dt_n e^{iP_j}
\\ &=&
\sum_{j=1}^n\frac{1}{2\pi}\int\hat{f}(t_1,\ldots,t_n)e^{i(t_1Q_1+\ldots+t_nQ_n)}it_jdt_1\ldots
dt_n e^{iP_j}
\\ &=&
\sum_{j=1}^n \frac{\partial f(Q_1,\ldots,Q_n) }{\partial
Q_j}e^{iP_j}de^{iP_j}
\\ &\cong&
\sum_{j=1}^n \frac{\partial f(Q_1,\ldots,Q_n) }{\partial Q_j}dQ_j
\end{eqnarray*}
which deforms to the classical first order differential calculus on
$\mathcal{S}(\mathbb{R}^n)$. We have thus showed the above result
for Schwartz class functions, combining the fact that
$\mathcal{S}(\mathbb{R}^n)$ is dense in $C_0(\mathbb{R}^n)$, we
actually have showed its relation with the classical calculus for
each $C_0(\mathbb{R}^n)$ functions.

\section {Quantum Heisenberg Algebra}
The quantum Heisenberg manifolds $D_{\mu,\nu}^c$, a continuous
field of $C^*$-algebra is first introduced by Rieffel {~\cite{R}},
recently, {~\cite{AE}} uses the Fell bundles method to get the
same quantum Heisenberg algebra as defined in {~\cite{R}}. The
construction here to get the $3$-dimensional quantum Heisenberg
manifolds follows from {~\cite{AE}}. Let
$$G=\left\{ \begin{pmatrix}  1 & y & z \\
                  0 & 1 & x \\
                  0 & 0 & 1  \end{pmatrix}:x,y,z\in\mathbb{R}\right\}$$
be the Heisenberg group of $3\times 3$ matrices. The Heisenberg
manifold $M_c$ is the quotient $G/H_c$ where $H_c$ is the discrete
subgroup of $G$ when $x, y, cz\in \mathbb{Z}$. To facilitate the
notation, we identify the above matrix as the vector $(x,y,z)\in
\mathbb{R}^{3}$, under this coordinate system, the multiplication
rule is simply
\begin{equation}
(x,y,z)(m,n,p)=(x+m,y+n,z+p+y m)
\end{equation}
 So $M_c$ can be described
as the quotient of the Euclidean space $\mathbb{R}^{3}$ by the
right action of $H_c$ given by above multiplication. Denote
$[x,y,z]$ as the quotient class of $(x,y,z)$ in $M_c$.
\[\phi_{(a,b)}[x,y,z]=[x+2b\mu,y+2b\nu,z+2b\nu x+2b^2\mu \nu+a/c]\]
defines an action of $\mathbb{R}^2$ on $M_c$. This action
generates the deformation data
$(\mathbb{Z},\gamma,\theta^\hbar)$(see definition 4.2
{~\cite{AE}}), the gauge action $\gamma$ is an action on the
circle group defined by
\[\gamma_{e^{2\pi i t}}([x,y,z])=[x,y,z+t/c]\] and the deformation
actions $\theta_\hbar$ of $\mathbb{Z}$ is defined by
\[\theta_{\hbar}([x,y,z])=[x+2\hbar \mu,y+2\hbar \nu,z+2\hbar \nu x+2 \hbar^2 \mu\nu].\]
For each $k\in\mathbb{Z}$, let
\begin{eqnarray*}
B_k & = & \{f\in C(M_c): \gamma_{e^{2\pi i \theta}}(f)=e^{2\pi i
\theta k}f, \forall \theta \in \mathbb{R} \}
\\ & = &
\{e^{2\pi i k cz}f(x,y): f(x+1,y)=e^{-2\pi i kc y}f(x,y), f(x,y)\in
C(\mathbb{R}\times \mathbb{T})\}.
\end{eqnarray*}
Take $f_k\in B_k$, $g_j\in B_j$, define the product and involution
by
\[f_k\times g_j=f_k\theta_{k\hbar}(g_j)\]
\[f_k^*=\theta^{-1}_{k\hbar}(\bar{f_k})\]
Then $C(M_c)_\gamma^\theta=(B_k,k\in\mathbb{Z},\times,*)$ becomes
a Fell bundle by keeping the linear, topological and norm
structure from $C(M_c)$. Theorem (8.3){~\cite{AE}} shows
$C(M_c)_\gamma^\theta\cong D_{\mu,\nu}^c$. We denote $N_{\hbar}$
as the von Neumann algebra of the weak operator closure of
$C(M_c)_\gamma^\theta$ on $L^2(\mathbb{R}\times \mathbb{T}^2)$.
For the time being, without further notice, let us assume
$\hbar=1$ or, what amounts to the same, that $\mu$ and $\nu$ are
replaced, respectively, by $\hbar\mu$ and $\hbar\nu$.

In order to study the structure of Heisenberg manifolds, one
possible way, as from the quantum tori, is to represent this $C^*$
algebra by generators. From (8) {~\cite{AE}}, we know that
$C(M_c)_\gamma^\theta$ is generated by $B_0,B_1$. As $B_0$ is
isomorphic to $C(\mathbb{T}^2)$, we get two generators $e^{2\pi
ix}$, $e^{2\pi i y}$, if having a third generator, the generator
$f_0(x,y)$ should not vanish at any point.

\newtheorem{hm}[guess]{Propostion}
\begin{hm}
Elements in $B_k(k\neq 0)$ always vanish somewhere.
\end{hm}
\begin{proof} We give the proof only for the case $k=1$, the other cases are
the same.

Otherwise, assume $0\neq f_0\in B_1$ and is also smooth, then
$e^{2\pi i c z}\frac{\partial f_0}{\partial x}\in B_1$, and
$f_0^{-1}\frac{\partial f_0}{\partial x} \in C(\mathbb{T}^2)$, hence
we can find a function $g\in C(\mathbb{T}^2)$ such that
$\frac{\partial f_0}{\partial x}=gf_0$. By solving this ODE, we have
$f_0(x,y)=C e^{\int_{-\infty}^x g(s,y)ds}$ where $C$ is some
constant. As $f_0\in B_1$, $f_0(x+1,y)=e^{-2\pi i c y}f_0(x,y)$,
this gives $e^{-2\pi i c y}=e^{\int_{0}^1 g(x,y)dx}$, but we don't
have such a function $g\in C(\mathbb{T}^2)$.

For an arbitrary $f\in B_1$, we find a smooth family $\{f_{n}\}_n\in
B_1$ which uniformly converges to $f$. For each $f_n$, it has a
vanishing point $x_n$. Since $f_n(x+1,y)=e^{-2\pi i cy}f_n(x,y)$,
there is no loss of generality in assuming $x_n\in [0,1]$. Further,
we can find a subsequence(we use the same $n$), such that
$x_n\rightarrow x^*$. Notice $f$ is uniformly continuous on $[0,1]$.
Given $\epsilon>0$, $\exists N, \delta$, such that when $n>N$ and
$|s-t|< \delta$,
\[|f_n(s)-f_n(t)|\leq |f_n(s)-f(s)|+|f(s)-f(t)|+|f(t)-f_n(t)|< \epsilon.\]
In particular,
$|f_n(x_m)-f_n(x^*)|<\epsilon$ when $|x_m-x^*|<\delta$ and $n>N$.\\
Let $n=m$ be large enough, we conclude $|f_n(x^*)|< \epsilon$. Let
$n\rightarrow \infty$, we have $|f(x^*)|<\epsilon$ which shows $x^*$
is a vanishing point of $f$.
\end{proof}

Above proposition tells us elements in $B_1$ always vanish
somewhere, we can't represent $C(M_c)_\gamma^\theta$ simply by
three generators. We have to find some other way to represent
$C(M_c)_\gamma^\theta$. Notice $C(M_c)_\gamma^\theta$ is the
closure of the direct sum $\bigoplus_{k\in \mathbb{Z}} B_k$, any
$f(x,y,z)\in C(M_c)_\gamma^\theta$ can be decomposed uniquely by
the Fourier transform on the $z$-coordinate, namely,
\[f(x,y,z)=\sum_{k\in \mathbb{Z}} e^{2\pi i kc z}f_k(x,y)\] where
$e^{2\pi i kc z}f_k(x,y)\in B_k$, and
\[ f_k(x,y)=\int_0^1 f(x,y,z)e^{-2\pi i k cz}dz\]
If we denote $\Phi(x,y,k)=f_k(x,y)$, it turns out $\Phi$ is a
function on $\mathbb{R}\times \mathbb{T}\times \mathbb{Z}$, which is
the original definition of $D_{\mu,\nu}^c$ by Rieffel {~\cite{R}}.

From {~\cite{R}} and above correspondence,
$\tau(f)=\int_0^1\int_0^1\int_0^1 f(x,y,z)dxdydz$ is a faithful
normal finite trace state and is invariant from the Heisenberg group
action. As $C(M_c)_\gamma^\theta$ is the closure of the direct sum
$\bigoplus_{k\in \mathbb{Z}} B_k$, we hence get
$L^2(C(M_c)_\gamma^\theta,\tau)=\bigoplus_{k\in
\mathbb{Z}}\bar{B}_k$, where $\bar{B}_k$ denotes the completion of
$B_k$ under $\tau$.

\newtheorem{iso}[guess]{Propostion}
\begin{iso} \label{bb}
$\bar{B}_k$ and $L^2(\mathbb{T}^2)$ is isomorphic.
\end{iso}
\begin{proof} Each element in $\bar{B}_k$ has the form $e^{2\pi
ikcz}f(x,y)$, where $f(x+1,y)=e^{-2\pi i kc y}f(x,y)$ and
$f(x,y+1)=f(x,y)$, which is easy to infer from the element from
$B_k$ and completion under $\tau$. So we define the linear map
$\tilde{U}: \bar{B}_k \rightarrow L^2(\mathbb{T}^2)$ by
\[\tilde{U} e^{2\pi ikcz}f(x,y)=e^{2\pi i c k x\{y\}}e^{2\pi ikcz}f(x,y).\]
Then,
\begin{eqnarray*}
& & \tau((e^{2\pi ikcz}f(x,y))^* \times e^{2\pi ikcz}f(x,y))
\\& =
& \tau(\overline{f(x-2k\mu,y-2k\nu)}f(x-2k\mu,y-2k\nu))
\\ & = & \int_{0}^1\int_{0}^1
\overline{f(x-2k\mu,y-2k\nu)}f(x-2k\mu,y-2k\nu)dx dy
\\& = &\int_{0}^1\int_{0}^1\overline{f(x,y)}f(x,y)dx dy.
\end{eqnarray*}
Hence $\langle e^{2\pi ikcz}f(x,y),e^{2\pi
ikcz}f(x,y)\rangle_{\bar{B}_k}=\langle \tilde{U}e^{2\pi
ikcz}f(x,y), \tilde{U}e^{2\pi
ikcz}f(x,y)\rangle_{L^2(\mathbb{T}^2)}$ and $\tilde{U}$ is an
isometry. As $\varphi_{m,n}=e^{2\pi i (mx+ny)}e^{2\pi i k c
z}e^{-2\pi i c k x \{y\}}\in \bar{B}_k$(see ~\cite{CSa}), and $\{
\tilde{U}\varphi_{m,n}\}_{m,n\in\mathbb{Z}}$ is an orthonormal
basis $L^2(\mathbb{T}^2)$, which implies $\tilde{U}$ is also
surjective.
\end{proof}
From above, we know $\{\varphi_{m,n}\}_{m,n \in \mathbb{Z}}$ is an
orthonomal basis of $\bar{B}_k$, but the basis doesn't derive the
multiplication structure from $\bar{B}_1$. Now let $U=e^{2\pi i
x},V=e^{2\pi i y}, W=e^{2\pi i cz}e^{-2\pi i c x\{y\}}$, where the
notation $\{x\}$ means the fraction part of $x$. Next proposition
will show $\{U^mV^nW^k\}_{m,n\in\mathbb{Z}}$ is an orthonomal
basis of $\bar{B}_k$.
\newtheorem{ba1}[guess]{Propostion}
\begin{ba1}
$\{U^mV^nW^k\}_{m,n\in\mathbb{Z}}$ is an orthonomal basis of
$\bar{B}_k$.
\end{ba1}
\begin{proof}Still denote $\varphi_{m,n}$ as in Proposition ~\ref{bb}, and let
\begin{eqnarray*}
\psi_{m,n} & = & U^mV^nW^k
\\ & = &
e^{2\pi i mx}e^{2\pi i ny}(\prod_{j=0}^{k-1}e^{2\pi i c(z+2j\nu
x+2j^2 \mu\nu)}e^{-2\pi i c(x+2j\mu)\{y+2j\nu\}})
\end{eqnarray*}
then $\varphi_{m,n}=\frac{\varphi_{m,n}}{\psi_{m,n}}\psi_{m,n}$,
and $|\frac{\varphi_{m,n}}{\psi_{m,n}}|=1$. $\psi_{m,n}$ is also
the orthonormal basis by the unitary transformation.
\end{proof}

\newtheorem{ba}[guess]{Corollary}
\begin{ba}
$\{U^mV^nW^k\}_{m,n,k\in \mathbb{Z}}$ is an orthonormal basis for
$L^2(C(M_c)_\gamma^\theta,\tau)$.
\end{ba}

\newtheorem{re}[guess]{Proposition}
\begin{re}
$N_{\hbar} \subset L^2(C(M_c)_\gamma^\theta,\tau)$.
\end{re}
\begin{proof}
This is because $\tau$ is a finite normal tracial state.
\end{proof}
As $U,V,W\in N_{\hbar}\subset L^2(C(M_c)_\gamma^\theta,\tau) $, we
conclude $N_{\hbar}$ is the weak operator closure generated by
$U,V,W$. The relation between $U,V,W$ is as follows:
\[ UV=VU, UW=e^{-4\pi i \mu}WU , VW=e^{-4\pi i \nu}WV,\]
\[UU^*=U^*U=VV^*=V^*V=WW^*=W^*W=1.\]
The above relation gives more, namely,  $N_{\hbar}$ is a 3
dimensional quantum tori. So we actually have proved the following
theorem.
\newtheorem{sam}[guess]{Theorem}
\begin{sam}
$C(M_c)_\gamma^\theta$(or $D_{\mu,\nu}^c$) is not a quantum tori,
while the von Neumann algebra $N_{\hbar}$ is  a 3-dimensional
quantum tori.
\end{sam}
In the language of the generators, another way to state the
faithful normal finite tracial state $\tau$ defined above is
\[\tau(\sum_{m,n,k} a_{m,n,k}U^mV^nW^k)=a_{0,0,0}.\]

In $D_{\mu,\nu}^c$, it has three canonical unbounded derivation
operators $D_1,D_2,D_3$ defined by
\[D_1 f=\frac{\partial f}{\partial x}=\sum_{k\in \mathbb{Z}} e^{2\pi i kc z}\frac{\partial f_k(x,y)}{\partial x},\]
\[D_2 f=\frac{\partial f}{\partial y}+x\frac{\partial f}{\partial z}=\sum_{k\in \mathbb{Z}} e^{2\pi i kc z}(\frac{\partial f_k(x,y)}{\partial y}+2\pi i kc x f_k(x,y)),\]
\[D_3 f=\frac{\partial f}{\partial z}=\sum_{k\in \mathbb{Z}} e^{2\pi i kc z} 2\pi i k cf_k(x,y).\]
We can now rephrase above unbounded derivation operators on
$N_{\hbar}$ in terms of generators as follows:
\[D_1 U= 2\pi i U, D_1 V=0, D_1 W=\sum_{k\neq 0}ck^{-1}V^k W,\]
\[D_2 U=0, D_2 V=2\pi i V, D_2 W=0,\]
\[D_3 U=0, D_3 V=0, D_3 W=2\pi i c W.\]
The summation above is in the weak operator convergence sense. The
elementary relation properties of these derivative operators are
listed below: \[[D_1,D_2]=D_3, [D_1,D_3]=[D_2,D_3]=0.\]

So far, we studied the structure of the quantum Heisenberg
manifolds, instead of using the quantum Heisenberg $C^*$-algebra
to apply the noncommutative differential calculus, we find it much
easier to use it in the von Neumann algebra case, as it is nothing
but a 3 dimensional quantum tori. The remainder of this section is
quite the same as from section ~\ref{tori} .

Suppose an element $a\in N_{\hbar}$ has the form
\[a=f(U,V,W)=\sum_{m,n,k}a_{m,n,k}U^mV^nW^k,\] then
\begin{eqnarray*}
& & [U^{k_1}V^{k_2}W^{k_3}, f(U,V,W)]
\\ & = &
\sum a_{m,n,k}(e^{4\pi i k_3(m\mu+n\nu)-4\pi i k(k_1
\mu+k_2\nu)}-1)U^mV^nW^k U^{k_1}V^{k_2}W^{k_3}
\\ & =&
(f(e^{4\pi i k_3\mu}U, e^{4\pi i k_3 \nu}V, e^{-4\pi
i(k_1\mu+k_2\nu)}W)-f(U,V,W))U^{k_1}V^{k_2}W^{k_3}
\end{eqnarray*}
The viewpoint of above derivation sheds some light on the changing
of the phase space in three directions, namely, with $4\pi i k_3
\mu$ in $x$-direction, $4\pi i k_3 \nu$ in $y$-direction, and $4\pi
i(k_1\mu+k_2 \nu)$ in $z$-direction. When $\hbar \rightarrow 0$, the
deformed derivation is as follows:
\[\lim_{\hbar\rightarrow 0}\frac{1}{\hbar}[W,f(U,V,W)]=(4\pi i \mu \frac{\partial f(U,V,W)}{\partial U}U +4\pi i \nu \frac{\partial f(U,V,W)}{\partial V}V)W,\]
\[\lim_{\hbar\rightarrow 0}\frac{1}{\hbar}[U,f(U,V,W)]=-4\pi i\mu \frac{\partial f(U,V,W)}{\partial W}WU,\]
\[\lim_{\hbar\rightarrow 0}\frac{1}{\hbar}[V,f(U,V,W)]=-4\pi i\nu \frac{\partial f(U,V,W)}{\partial W}WV.\]

\bibliographystyle{amsplain}

\end{document}